\documentclass{amsart}
\usepackage{amssymb,amsthm, amsmath, amscd, mathrsfs}
\newcommand{\nc}{\newcommand}

\nc{\KLlam}{\mathcal{M}^{\slam}_0}
\nc{\Catlam}{\RZlam\grmod}
\nc{\RZlam}{R_{\Center}^{\slam}}
\nc{\chilam}{\chi_{\slam}}
\nc{\Vg}{\V_{\g,\cri}}
\nc{\kc}{\m_{\chi_{\lam}}}
\nc{\Clam}{\mathcal{M}_{\cri}^{\slam}}
\nc{\Cp}[1]{\mathcal{M}_{\cri}^{\lam,#1}}
\nc{\Catp}{\mathcal{M}_{\cri}}
\nc{\Catc}{{\mathcal{M}_{\cri}^0}}
\nc{\RZ}{R_{\Center}}\nc{\VV}{\mathbb{W}}
\nc{\UZ}{R_{\Center}^-}
\nc{\ZZ}{\mathfrak{Z}(\V_{\g,\cri})}
\nc{\CZ}{\mathcal{C}_{\mathcal{Z}}}
\nc{\Lcalflys}{\mathcal{L}}
\nc{\bbK}{\mathbb{K}}
\nc{\bbM}{\mathbb{M}}
\nc{\OO}{\mathbb{O}}
\nc{\Cen}{\Center_{\g,-h\che}}
\nc{\bnp}{\bn_+}
\nc{\bnm}{\bn_-}
\nc{\snp}{\sn_+}
\nc{\snm}{\sn_-}
\nc{\OpG}{\Op_{{}^L G}}
\nc{\calR}{\mathcal{R}}
\nc{\Ope}{\Op_{{}^L G}(D^{\times})}
\nc{\twoheadra}{\twoheadrightarrow}
\nc{\calM}{\mathcal{M}}
\nc{\calK}{\mathcal{K}}
\nc{\sups}{\supset}
\nc{\sub}{\subset}
\nc{\prroots}{\roots_+^{\re}}
\nc{\wchi}{\tilde{\chi}}
\nc{\np}[1]{\n_{#1}^+}
\nc{\nm}[1]{\n_{#1}^-}
\nc{\HH}{H_{\BRST}}
\nc{\Paf}[1]{P^{+}_{#1}}
\nc{\biglam}{\widehat{\lam}}
\nc{\G}{\mathcal{G}}
\nc{\germ}[1]{\mathfrak{#1}}
\nc{\cprime}{$'$}
\nc{\Fi}{\mathbb{F}}
\nc{\Sy}{\mathfrak{S}}
\nc{\gl}{\mathfrak{gl}}
\nc{\Wfin}[1]{\W^{\fin}_{#1}(\g,\chi)}
\nc{\hookra}{\hookrightarrow}
\nc{\ban}[1]{\rm{(}#1\rm{)}}
\nc{\longra}{\longrightarrow}
\nc{\Um}{\U^-}
\nc{\CClam}{C}
\nc{\zero}{0}
\nc{\I}{\mathsf{I}}
\nc{\nmap}{\Psi}
\nc{\WW}{\mathsf{W}}
\nc{\vQ}{\mathsf{Q}}
\nc{\vchi}{\mathsf{\chi}}
\nc{\sig}{\sigma}
\nc{\TT}{T}
\nc{\sE}{\bar{E}}
\nc{\Lone}{T^*}
\nc{\Lnew}{T_1}
\nc{\Lold}{T_1^{\old}}
\nc{\nno}{\nonumber}
\nc{\omu}{\bar{v}_{\mu}}
\nc{\gU}{\a_{U''}}
\nc{\olam}{\overline{|\lam\ket}}
\nc{\bCmu}{K(\mu)}
\nc{\bClam}{K(\lam)}
\nc{\ttp}{\tp}
\nc{\FF}{\mathsf{F}}
\nc{\nW}{\bar{\W}_+}
\nc{\newWU}{\wt{U_k(\bg')\* \bCl}_{\new}}
\nc{\oldWU}{\wt{U_k(\bg')\* \bCl}_{\old}}
\nc{\newU}{U_k(\bg')\* \bCl_{\new}}
\nc{\ra}{\rightarrow}
\nc{\Hp}[1]{H^{#1}}
\nc{\UV}{\U(\V)}
\nc{\wJ}{\widehat{J}}
\nc{\wBB}{\wt{\B}}
\nc{\UN}[1]{\U_N(#1)}
\nc{\QO}{Q}
\nc{\AN}[1]{\mathscr{A}_N(#1)}
\nc{\adMod}{\text{-}\mathsf{Mod^{ad}}}
\nc{\GrMod}{\text{-}\mathsf{Mod^{ad}}}
\nc{\Levi}[1]{L_{#1}(\sg,\chi)}
\nc{\bCgoo}[1]{B_{#1}(\sg,\chi)''}
\nc{\bCgo}{\bCg''}
\nc{\bCGo}{\bCg''}
\nc{\wUU}{\wt{\UU}}
\nc{\QN}{\U(\V)/\wt{\II}_N}
\nc{\QNS}[1]{\U(\V_{#1})/\wt{\II}_N^{(#1)}}
\nc{\LV}{L(\V)}
\nc{\wh}{\widehat}
\nc{\adQ}{\ad Q }
\nc{\dV}{Q_{(0)}}
\nc{\wdV}{\wt{\QO}}
\nc{\Basis}{\mathcal{B}}
\nc{\UU}{\mathbf{U}}
\nc{\bU}{\mathbb{U}}
\nc{\II}{\mathscr{I}}
\nc{\wt}{\widetilde}
\nc{\Fin}{\mathsf{Fin}}
\nc{\Vkg}{V_k(\sg)}
\nc{\parw}{\partial_w}
\nc{\wU}{\widetilde{U}}
\nc{\Wg}[1]{\W_{#1}(\sg)}
\nc{\wc}{\widetilde{c}}
\nc{\twlam}{\overline{t_{-\srho\che}\circ \lam}}
\nc{\bge}{\bg^e}
\nc{\bgf}{\bg^f}
\nc{\sge}{\sg^e}
\nc{\sgf}{\sg^f}
\nc{\WN}{\mathbf{N}}
\nc{\bt}{\bh}
\nc{\bCl}{\Cl}
\nc{\bCg}{C_k(\sg)}
\nc{\sCg}[1]{C(\g_{#1})}
\nc{\sCgo}[1]{C(\g_{#1})''}
\nc{\sN}{\bar{N}}
\nc{\sG}{\bar{G}}
\nc{\smu}{\bar{\mu}}
\nc{\salpha}{\bar{\alpha}}
\nc{\hatJ}{\widehat{J}}
\nc{\sA}{\bar{A}}
\nc{\Cat}{{\mathcal{C}}}
\nc{\Lie}[1]{{\mathfrak{L}(#1)}}
\nc{\sw}{\bar{w}}
\nc{\sC}{\bar{C}}
\nc{\dBGG}{\dot{\BGG}}
\nc{\tri}{\triangle}
\nc{\GG}{\mathsf{Y}}
\nc{\Prnongen}{{\textsl{Pr}}_{\kappa,\mathrm{nondeg}}}
\nc{\Center}{{\mathcal{Z}}}
\nc{\Wcat}{\Wg\text{-}\mathsf{grmod}}
\nc{\WCat}{\Wcat}
\nc{\fmap}{\rightsquigarrow }
\nc{\fdomain}{{\mathcal{C}}}
\nc{\D}{{\mathsf{D}}}

\nc{\whJ}{\widehat{J}}
\nc{\bd}{Q}
\nc{\blam}{\lam}
\nc{\dchi}{\chi}
\nc{\sV}{\bar{V}}
\nc{\sCl}{{\Cl}}
\nc{\sH}{\bar{H}}
\nc{\schi}{\bar{\chi}}
\nc{\sL}{\bar{L}}
\nc{\sM}{\bar{M}}
\nc{\sd}{\bar{Q}}
\nc{\slam}{\bar{\lam}}

\newcommand{\V}{{\mathbb V}}

\newcommand{\M}{{\mathcal{M}}}
\newcommand{\B}{{\mathbf{B}}}

\newcommand{\N}{{\mathbb{N}}}
\newcommand{\U}{{\mathfrak{U}}}

\newcommand{\isomap}{\overset{\sim}{\rightarrow} }
\newcommand{\W}{{\mathscr{W}}}

\newcommand{\sroots}{\bar{\Delta}}
\newcommand{\roots}{\Delta}

\newcommand{\rroots}{\Delta^{\rm re}}

\newcommand{\m}{{\mathfrak{m}}}

\renewcommand{\a}{{\mathfrak{a}}}

\newcommand{\Cl}{{\mathcal{C}l}}

\newcommand{\KL}{\BGG^{\textKL}_{\cri}}

\newcommand{\C}{{\mathbb C}}
\newcommand{\Z}{{\mathbb Z}}

\newcommand{\inv}{^{-1}}
\newcommand{\dual}[1]{{#1}^*}
\newcommand{\lam}{\lambda}
\newcommand{\Lam}{\Lambda}

\renewcommand{\*}{{\otimes}}
\newcommand{\+}{\mathop{\oplus}}
\newcommand{\h}{ {\mathfrak h}}
\newcommand{\g}{ {\mathfrak g}}

\newcommand{\che}{^{\vee}}
\newcommand{\bra}{{\langle}}
\newcommand{\ket}{{\rangle}}

\DeclareMathOperator{\grmod}{-grmod}
\DeclareMathOperator{\Spec}{Spec}

\DeclareMathOperator{\textKL}{KL}
\DeclareMathOperator{\BRST}{BRST}
\DeclareMathOperator{\cri}{crit}
\DeclareMathOperator{\Op}{Op}

\DeclareMathOperator{\Zhu}{{Zh}}

\DeclareMathOperator{\tp}{top}

\DeclareMathOperator{\End}{End}

\DeclareMathOperator{\Hom}{Hom}

\DeclareMathOperator{\ad}{ad}

\DeclareMathOperator{\Ann}{Ann}
\DeclareMathOperator{\fin}{fin}

\DeclareMathOperator{\ch}{ch}

\DeclareMathOperator{\rank}{rank}

\DeclareMathOperator{\new}{new}
\DeclareMathOperator{\old}{old}

\newcommand{\sg}{ \bar{\mathfrak g}}
\newcommand{\sh}{\bar{ \h}}
\newcommand{\sn}{\bar{\mathfrak{n}}}

\newcommand{\sP}{\bar P}
\newcommand{\sW}{\bar{W}}
\newcommand{\srho}{ \bar{\rho}}

\newcommand{\sproots}{\bar{\Delta}_+}

\newcommand{\bg}{{\g}}
\newcommand{\bh}{ \h}
\newcommand{\bn}{\n}

\DeclareMathOperator{\re}{re}

\renewcommand{\bigg}{{\mathfrak g}}

\newcommand{\n}{{{\mathfrak{n}}}}

\newcommand{\BGG}{{\mathcal O}}

\title[Affine Kac-Moody   Lie algebras 
at the critical level]
{Characters  of  representations of affine Kac-Moody Lie algebras 
at the critical level
}
\author{Tomoyuki  Arakawa}
\address{Department of Mathematics, 
Nara Women's University, Nara 630-8506, JAPAN}
\email{arakawa@cc.nara-wu.ac.jp}
\thanks{The author is partially supported 
by the JSPS Grant-in-Aid for Young Scientists (B)
No.\ 17740006}
\begin{document}
\theoremstyle{plain}
\newtheorem{Th}{Theorem}
\newtheorem{MainTh}{Main Theorem}
\newtheorem{Pro}[Th]{Proposition}
\newtheorem{ProDef}[Th]{Proposition and Definition}
\newtheorem{DefPro}[Th]{Definition and Proposition}
\newtheorem{Lem}[Th]{Lemma}
\newtheorem{Co}[Th]{Corollary}
\newtheorem{Facts}[Th]{Facts}
\theoremstyle{definition}
\theoremstyle{remark}
\newtheorem{Def}[Th]{Definition}
\newtheorem{Rem}[Th]{Remark}
\newtheorem{Conj}{Conjecture}
\newtheorem{Claim}{Claim}
\newtheorem*{ClaimNN}{Claim}
\newtheorem{Notation}{Notation}
\newtheorem{Ex}[Th]{Example}
\newcommand{\st}{{\mathrm{st}}}
\maketitle
\section{Introduction and Main Results}
\subsection{}Let $\sg$ be a complex simple Lie algebra of rank $l$,
$\bigg$ non-twisted affine Kac-Moody Lie algebra associated with
$\sg$:
\begin{align}
 \bigg=\sg\* \C[t,t\inv]\+\C K\+ \C D.
\end{align}
The commutation relations of $\bigg$
are given by
the following.
\begin{align*}
& [X(m),Y(n)]=[X,Y](m+n)+m\delta_{m+n,0}(X|Y)K,\\
&[D,X(m)]=mX(m),
\quad [K, \bg]=0
\end{align*}
for 
$X, Y\in \sg$,
$m,n\in \Z$,
where 
$X(m)=X\* t^m$
with $X\in \sg$ and $m\in \Z$
and 
$(\cdot |\cdot)$ is the normalized invariant inner product of $\sg$.
We identify $\sg$ with $\sg \* \C\subset \bg$.
Fix the triangular decomposition $\sg=\snm\+ \sh\+ \snp$,
and  the Cartan subalgebra of $\bg$
as  $\bh=\sh\+ \C K \+ \C D$.
We have $\dual{\bh}=\dual{\sh}\+ \C \Lam_0\+ \C \delta$,
where $\Lam_0$ and $\delta$ are elements dual to
$K$ and $D$,
respectively.

Let $L(\lam)$ be the irreducible highest weight 
representation of $\g$ of highest weight
$\lam\in \dual{\bh}$ with respect to the standard triangular
decomposition
$\bg=\bnm\+ \bh\+ \bnp$,
where 
\begin{align*}
\bnm=\snm \+ \sg\* \C[t\inv] t\inv,\quad
\bnp=\snp \+\sg \* \C[t] t.
\end{align*}
The central element $K$
acts on $L(\lam)$ as  the multiplication
by the constant $\bra \lam,K\ket$,
which is called the {\em level} of $L(\lam)$.
The level $\bra \lam,K\ket =-h\che$ 
is called {\em critical},
where $h\che$
is 
the dual  Coxeter number of $\sg$.

\subsection{}
Let $\ch L(\lam)$
be the formal character of $L(\lam)$:
\begin{align*}
\ch L(\lam)=\sum_{\mu\in \dual{\bh}} e^{\mu}\dim_{\C}L(\lam)^{\mu},
\end{align*}
where $L(\lam)^{\mu}$ is the weight space of $L(\lam)$
of weight $\mu$.

 The Weyl-Kac  character formula \cite{Kac74}
gives an
explicit formula of $\ch L(\lam)$
in the case that
$L(\lam)$ is
an  integrable representations of $\g$.
It is known that
Kac-Wakimoto admissible representations \cite{KacWak88,KacWak89}
also have
the
 Weyl-Kac  type character formulas.
The celebrated 
{Kazhdan-Lusztig conjecture} \cite{KazLus79}
(proved by \cite{BeuBer81,BryKas81})
has been generalized to $\bg$
by Kashiwara-Tanisaki \cite{KasTan95,KasTan96,KasTan98,KasTan00}
and Casian \cite{Cas96}.
As a result
the character $\ch L(\lam)$ is known for 
{\em any}
$L(\lam)$
provided that
its level
is
{\em not}
critical (see \cite{KasTan00} for the most general formula).

\subsection{}
On the contrary 
not much  is known 
about 
the characters of $L(\lam)$
at the critical level.
It seems that
 only the
{\em generic case} is known,
that is,
the case that
$\lam$ satisfies the condition that
\begin{align}
 \bra \lam+\rho,\alpha\che\ket \not \in \N
\quad \text{for all $\alpha\in \rroots_+$},
\label{eq:2007-06-03-23-34}
\end{align}
where  $\rroots_+$ is the set of positive real roots of $\bg$,
 $\rho=\srho+h\che \Lam_0$,
 $\srho=1/2\sum_{\salpha\in \sroots_+}\alpha$
and $\sroots_+\subset \rroots_+$ is the set of positive roots of $\sg$.
In this case 
the Kac-Kazhdan conjecture \cite{KacKaz79}
(which is a theorem proved by 
\cite{Hay88,GooWal89,FeuFre88,Ku89}\footnote{See also
\cite{Mat03,Fre05,Ara06}.})
gives the following character formula
of $L(\lam)$:
\begin{align}
 \ch L(\lam)=\frac{e^{\lam}}{\prod\limits_{\alpha\in
 \rroots_+}(1-e^{-\alpha})}.
 \label{eq:2007-06-03-23-43}
\end{align}

By  the existence of the Wakimoto modules at the 
critical level \cite{Wak86,FeuFre90,Fre05},
it follows that
the irreducible representation $L(\lam)$ at the critical level
in general 
has a character equal  to or smaller 
 than the right hand side of
\eqref{eq:2007-06-03-23-43}.

\subsection{}
In this paper 
%
 we study the irreducible highest weight representations of $\bg$
at the critical level
which are integrable over $\sg$.
Denote by  $\slam$  the restriction  of $\lam\in \dual{\h}$ 
to 
$\sh$.
Set \begin{align*}
     \Paf{\cri}=\{\lam\in\dual{\bh};
\slam\in \sP^+, \ \bra \lam,K\ket=-h\che
\},
    \end{align*}
where $\sP^+$ the set of integral dominant weights of $\sg$:
\begin{align}
\sP^+=\{\slam\in \dual{\sh};
\bra \slam,\alpha\che\ket\in \Z_{\geq 0}\text{ for all }
\alpha\in \sproots\}.
\end{align}
The $L(\lam)$ at the critical
level
is integrable over $\sg$
if and only if $\lam$ belongs to $\Paf{\cri}$.

We have the following result.
\begin{Th}\label{Th:Main}
Let $\lam\in \Paf{\cri}$.
The character of $L(\lam)$
is given by
\begin{align*}
 \ch L(\lam)=
\displaystyle{\frac{\sum\limits_{w\in \sW}(-1)^{\ell(w)}e^{w\circ \lam}}{
\prod\limits_{\alpha\in \sroots_+}(1-q^{-\bra \lam+\rho,\alpha\che\ket})
\prod\limits_{\alpha\in \rroots_+}(1-e^{-\alpha})}},
\end{align*}
where $q=e^{\delta}$,
$\sW$ is the Weyl group of $\sg$,
$\ell(w)$ is the length of $w$
and $w\circ \lam=w(\lam+\rho)-\rho$.
\end{Th}
\subsection{}
Recently
the representations of
$\bg$ at the critical level have been
studied in detail by Frenkel and Gaitsgory
\cite{FreGai04,FreGai05b,FreGai06,FreGai07}
in the view point of the 
{\em geometric Langlands program}.
Our original motivation was to
 confirm\footnote{After finishing 
 this paper we were notified that
Frenkel and Gaitsgory have
recently proved Conjecture 5 of \cite{Fre06}
and that Theorem  \ref{Th:Main} was  known to E. Frenkel.
}
Conjecture 5 of \cite{Fre06}
in the case that the 
{\em opers} are  ``graded''
 (see Theorem \ref{Th:chonj-tru})
by applying the method of 
the {\em quantum Drinfeld-Sokolov reduction}
\cite{FF90,FeiFre92,FreBen04} 
(cf.\ \cite{FKW92,KacRoaWak03,KacWak04,Ara04,Ara05,Ara06,Ara07}).
Theorem \ref{Th:Main} has been obtained as a byproduct of the proof.

\section{Endmorphism rings,
duality and the character  formula}
\subsection{The category $\KL$}
Denote by $\KL$ the full subcategory of 
the category of $\bg$-modules 
consisting of objects $M$ such that the following hold:
(1) $K$ acts on $M$ as the multiplication by $-h\che$,
(2) $D$ acts on $M$ semisimply:
$M=\bigoplus_{d\in \C}M_d$,
where $M_d=\{m\in M; D m=d m\}$,
(3) $\dim M_d<\infty $ for all $d$,
(4) there exists a finite subset $\{d_1,\dots ,d_r\}$ 
of $\C$ such that 
$M_d=0$ unless $d\in \bigcup_{i=1^r}d_i-\Z_{\geq 0}$.

Any object $M$ of $\KL$
admits a weight space decomposition:
$M=\bigoplus_{\lam}M^{\lam}$,
$M^{\lam}=\{m\in M; h m=\bra \lam,h\ket m\, \forall h\in \h\}$.
We set $\ch M=\sum_{\lam\in \dual{\bh}}e^{\lam}\dim M^{\lam}$.

The {\em Weyl module}
\begin{align}
V(\lam)=U(\bg)\*_{U(\sg\* \C[t]\+ \C K)} E(\lam)
\end{align}
with $\lam\in \Paf{\cri}$
belongs to $\KL$.
Here
$E(\lam)$ is the irreducible finite-dimensional
representation of $\sg$
of highest weight $\slam$,
considered as a $\sg\* \C[t]\+ \C K \+ \C D$-module
on which $\sg\* \C[t]t$ acts by zero,
and
 $K$ and $D$ act  as the multiplication
by $-h\che$ and $\bra \lam,D\ket$,
respectively.
By the Weyl character formula one has
\begin{align}\label{eq:2007-05-28-00-08}
 \ch V(\lam)=
\displaystyle{\frac{\sum_{w\in \sW}(-1)^{\ell(w)}e^{w\circ \lam}}{
\prod_{j\geq 1}(1-q^{-j})^{\ell }
\prod_{\alpha\in \rroots_+}(1-e^{-\alpha})}}.
\end{align}
The $V(\lam)$ has $L(\lam)$
as its unique  simple quotient.

\subsection{The derived algebra $\bg'$ of $\bg$}
Let $\g'$ be the derived algebra of  $\bg$:
\begin{align}
 \bg':=[\bg,\bg]=\sg\*\C[t,t\inv]\+ \C K.
\end{align}
One sees that each $L(\lam)$ remains irreducible over $\bg'$.
\subsection{The vertex algebra associated with $\sg$ at the critical level}
The vacuum Weyl module   
\begin{align}
 \V_{\g,\cri}:=V(-h\che \Lam_0) 
\end{align}
has
the natural structure of
vertex algebras (see eg.\ \cite{Kac98,FreBen04}).
The $\V_{\g,\cri}$ is called the
{\em  universal
affine 
vertex algebra 
associated with $\sg$ at the critical level}.
Each object of $\KL$ can be regarded as a
$\V_{\g,\cri}$-module.

For a vertex algebra $V$ in general
we denote by  
\begin{align}
Y(a,z)=\sum_{n\in \Z}a_{(n)}z^{-n-1}\in (\End V)[[z,z\inv]]
\end{align}
 the quantum field corresponding to $a\in V$.
Also,
we write $\Zhu(V)$ for the {\em Zhu algebra} \cite{FreZhu92,Zhu96}
(see also \cite{NagTsu05}) of $V$.
One knows by  \cite{FreZhu92} 
that there is a natural isomorphism
\begin{align}
 \Zhu(\V_{\g,\cri})\cong U(\sg).
\end{align}

 \subsection{Feigin-Frenkel's theorem}
The vertex algebra 
$\Vg$ has a large center \cite{Hay88,FeiFre92},
which we denote by $\ZZ$:
\begin{align*}
 \ZZ&=
\{a\in \V_{\g,\cri};
a_{(n)}v=0\text{ for all }n\in \Z_{\geq 0},\ v\in \V_{\g,\cri}\}.
\end{align*}
Then 
\begin{align*}
\ZZ
=
\{a\in \V_{\g,\cri};
[a_{(m)}, v_{(n)}]=0\text{ for all }m, n\in \Z,\ v\in \V_{\g,\cri}\}.
\end{align*}
Let $a\in \ZZ$
and $n\in \Z$.
The action of $a_{(n)}$ 
on $M\in \KL$ commutes with the action of $\bg'$.
Therefore each $a_{(n)}$ 
acts as the multiplication by a constant on $L(\lam)$.

Let $I=\{1,2,\dots, l\}$ ($l=\rank \sg$),
$\{d_i; i\in I\}$ the exponents of $\sg$,
 $\Center(\sg)$ the center of the universal enveloping 
algebra $U(\sg)$ of $\sg$.

There is a remarkable realization of 
$\ZZ$
due to Feigin and Frenkel \cite{FeiFre92}
as a chiralization of  Kostant's  {\em Whittaker model}
 \cite{Kos78} of  $\Center(\sg)$
(for  details, see \cite{FreBen04}).
As a result
one has the following description of $\ZZ.$
  
\begin{Th}[E. Frenkel and B. Feigin \cite{FeiFre92}]\label{Th:FF}
There exist  homogeneous vectors 
$p^i\in \ZZ\cap (\Vg)_{-d_i-1}$
with $i\in I$
that generate
a
PBW basis of $\ZZ$:
that is,
 there is a linear isomorphism
\begin{align*}
\begin{array}{ccc}
  \C[p^i_{(-n)}; i\in I, n\in \Z_{\geq 1}]
&
\isomap & \ZZ.
 \\
a&\mapsto & a|0\ket,
\end{array}\end{align*}
where $|0\ket $ is the highest weight vector of $\Vg$.
\end{Th}

\subsection{The linkage  principle for $\KL$ }
Let $p^i
$ with $i\in I$ be 
generators of  $\ZZ$ as in Theorem 
\ref{Th:FF}.
(We have $p^i=p^i_{(-1)}|0\ket$.)
Write
\begin{align}
 Y(p^i,z)=\sum_{n\in \Z}p^i_n z^{-n-d_i-1},
\end{align}
so that
\begin{align}
 [D, p^i_n]=n p^i_n
\label{eq:2007-05-29-14-44}
\end{align}
on $M\in \KL$.
Set
\begin{align}
 \RZ=\C[p^i_n; i\in I,\ n\in \Z].
\end{align}
An object $M$ of $\KL$ 
is regarded as a $\RZ$-module naturally.

There is a natural map
\begin{align}
   \C[p^i_0;i\in I] \ra \Zhu(\ZZ),
\end{align}
which is actually an isomorphism.
To be precise  we have the following assertion.
\begin{Th}\label{Th:iso-center}
The natural map $\Zhu(\ZZ)\ra \Zhu(\V_{\g,\cri})=U(\sg)$
is injective and its image coincides with
$\Center(\sg)$.
\end{Th}
See \cite{Ara07} for a proof of
Theorem \ref{Th:iso-center}.  
We shall identify $\Zhu(\ZZ)$ with $\Center(\sg)\subset U(\sg)$
through Theorem \ref{Th:iso-center}.

Let $o(p^i)$ be the image of $p^i_0$ in
$\Zhu(\ZZ)$.
If $v_{\lam}\in M\in \KL$
is annihilated by
$\bnp$, 
then one has
\begin{align}
 p^i_0 v=o(p^i)v=\bar\chi_{\slam}(o(p^i))v,
\end{align}
where 
\begin{align*}
\bar\chi_{\slam}: \Center(\sg)\ra \C
\end{align*}
 is the evaluation
of $\Center(\sg)$ at the Verma module of $\sg$ of highest weight
$\slam$.

Because $p_0^i$ commutes with the action of 
$\bg'$ 
on $M\in \KL$,
the following assertion follows
immediately.
\begin{Pro}\label{Pro:2007-05-22-22-22-24}
 If $L(\mu)$ appears in the local composition factor
of $V(\lam)$ then $\mu=\lam-n\delta$ for some $n\in \Z_{\geq 0}$.
\end{Pro} 
 
\subsection{The conjecture
of Frenkel and Gaitsgory for graded opers}
Let $\slam\in \sP^+$.
The character $\schi_{\slam}$
 naturally extends to
the graded central character 
of $\ZZ$,
that is,
to
the ring homomorphism
\begin{align}
\chilam : \RZ\ra \C
\end{align}
defined by 
\begin{align}
 \chilam (p^i_n)=\begin{cases}
		       \bar\chi_{\slam}(o(p^i_0))&\text{if }n=0,\\
0&\text{if }n\ne 0.
		      \end{cases}
\end{align}
The $\ker \chilam\cdot V(\lam)$
is a submodule of $V(\lam)$.
One has
\begin{align}
\ker\chilam\cdot V(\lam)=\sum_{n>0,\ i\in I}U(\bnm)p_{-n}^i |\lam\ket,
\label{eq:2007-05-23-003}
\end{align}
where 
 $|\lam\ket$ is  the highest weight vector of $V(\lam)$.
Thus $\ker \chilam\cdot V(\lam)$ is a proper submodule of $V(\lam)$
which lies in $\KL$.
Hence  there is a following exact sequence in $\KL$:
\begin{align}
V(\lam)/\ker \chilam\cdot V(\lam)\ra L(\lam)\ra 0.
\label{eq:2007-05-23-07-31}
\end{align}

The following assertion is clear.
\begin{Pro}
Any vector of $ L(\lam)$ is annihilated by $\ker \chilam$.
\end{Pro}
Let $\KLlam$ be the full subcategory 
$\KL$ consisting of objects $M$
such that $\ker \chilam\cdot M=0$.
Any simple object of $\KLlam$ 
is isomorphic to  $L(\mu)$ with
$\mu\in\Paf{\cri}$ such that $\bar \mu=\slam$
(thus all the simple modules are 
mutually isomorphic as $\g'$-modules).

The following striking assertion 
was conjectured by 
Frenkel and Gaitsgory (announced in \cite{Fre06}).
\begin{Conj}[E. Frenkel and D. Gaitsgory]\label{Conj-FG}$ 
$
\begin{enumerate}
 \item 
The category $\KLlam$ is semisimple
for any $\slam\in \sP^+$.
\item For each $\lam\in \Paf{\cri}$
there is
 an isomorphism
 $V(\lam)/\ker \chilam\cdot V(\lam)\cong L(\lam)$.
\end{enumerate}
\end{Conj}
\begin{Rem}\label{Rem:originial-conj}
$ $

\begin{enumerate}
\item 
By
the ``Langlands duality'' 
\cite{FeiFre92},
the $\chilam$
can be considered as an
element of
the {\em ${}^L G$-oper}
$\OpG(D^{\times})$ \cite{BeiDri04}
on the punctured disk $D^{\times}$,
which is ``graded''.
The original conjecture 
(Conjecture 5 of \cite{Fre06})
of Frenkel and Gaitsgory
is more general and applies to
any (not necessarily graded)
central character
$\chi$
(i.e.\ to any element of $\Op_{{}^L G}^{\lam}$,
see Remark \ref{Rem:reg-op} and 
\cite{Fre06}).
 \item In the case that $\slam=0$,
Conjecture \ref{Conj-FG}
follows from 
Theorem 6.3 of  \cite{FreGai04}
(applied to the graded oper $\chi_{0}$).
\end{enumerate} 
\end{Rem}

\subsection{Endmorphism rings of Weyl modules}
Let $\lam\in \Paf{\cri}$.
Recall that $|\lam \ket $ denotes 
the highest weight vector of $V(\lam)$.
Define
\begin{align}
 \RZlam=\RZ/\Ann_{\RZ}|\lam\ket.
\end{align}
Note that $\RZlam$ is naturally graded by $D$:
\begin{align}
 \RZlam=\bigoplus_{d\in -\Z_{\geq 0}}
(\RZlam)_{d},
\quad
(\RZlam)_{d}=\{a\in \RZlam; [D,a]=d a\}.
\end{align}

There is a 
natural algebra homomorphism 
\begin{align}
 \RZlam\ra \End_{U(\g')}(V(\lam)).
\label{eq:2007-06-04-00-52}
\end{align}
If $\slam =0$,
then $\RZ^0\cong \ZZ$.
In this case it is known by \cite{FeiFre92,Fre05}
that \eqref{eq:2007-06-04-00-52}
gives an isomorphism
\begin{align}
 \ZZ\cong \End_{U(\g')}(\Vg).
\end{align}
This is true for any $\lam\in \Paf{\cri}$.
\begin{Th}\label{Th:end-ring}
Let $\lam\in \Paf{\cri}$.
 \begin{enumerate}
  \item \label{item:endo-iso}
The map \eqref{eq:2007-06-04-00-52} 
gives the isomorphism 
$\RZlam\isomap \End_{U(\g')}(V(\lam))$.
\item \label{item:ch-of-end}
Set $\ch \RZlam=\bigoplus_{d\in \C}q^d
\dim (\RZlam)_d$.
Then 
\begin{align*}
 \ch \RZlam
=\displaystyle
{\frac{
\prod_{\alpha\in \sroots_+}(1-q^{-\bra \lam+\rho,\alpha\che\ket })}
{\prod_{j\geq 1}(1-q^{-j})^{\ell}}.
}
\end{align*}
 \end{enumerate}
\end{Th}
Theorem \ref{Th:end-ring} 
was obtained earlier in \cite{FreGai07}.
In \cite{Ara08} we give an
independent 
proof of Theorem \ref{Th:end-ring}
by  the method of 
quantum Drinfeld-Sokolov reduction.
\begin{Rem}\label{Rem:reg-op}
According to 
Frenkel and Gaitsgory \cite{FreGai07},
one has
\begin{align*}
\Spec \RZlam\cong \Op_{{}^L G}^{\slam}.
\end{align*}
where $\Op_{{}^L G}^{\slam}$ is a certain sub-pro-variety of
$\OpG(D^{\times})$ described in  \cite{Fre03,FreGai06}
(cf.\ (7.17) of \cite{Fre06}).
\end{Rem}
\subsection{An equivalence of  categories}
Let $\M^{\slam}$ be the full subcategory
of $\KL$ consisting of 
objects $M$ that are annihilated by 
\begin{align*}
 p_n^i-\chilam(p_n^i)
\quad \text{for all }i\in I,\ n\geq 0.
\end{align*}
Then $V(\lam)$, $L(\lam)\in \M^{\slam}$.
Also,
$\KLlam$ is a full subcategory of $\M^{\slam}$.

Let $\Catlam$ be the 
full subcategory of the category of graded
$\RZlam$-modules 
consisting of objects
$X=\bigoplus_{d\in \C}X_d$ 
($(\RZlam)_d\cdot X_{d'}\subset X_{d+d'}$)
such that 
(1) $\dim X_d<\infty$ for all $d\in \C$;
(2) there exists a finite subset $d_1,\dots,d_r\subset \C$
such that $X_d=0$ unless
$d\in \bigcup d_i-\Z_{\geq 0}$.
Then
any simple object of $\Catlam$ is isomorphic  to
\begin{align}
\C_{\chilam}:=\RZlam/\ker\chilam \cdot \RZlam
\end{align}
as $\RZlam$-modules.

Set
\begin{align}
& F(M)=\Hom_{U(\g')}(V(\lam),M)
\end{align}
for $M\in \M^{\slam}$.
Then 
$F(M)\cong M^{\bnp}$.
The  $F(M)$ is naturally a graded $\RZlam$-module: 
\begin{align}
F(M)=\bigoplus_{d\in \Z}F(M)_d
\end{align}
where $F(M)_d= M^{\bnp}\cap M_d$.
Thus
$F$ defines a functor from $\M^{\slam}$ to $\Catlam$.

Next let
\begin{align}
G(X)=V(\lam)\*_{\RZlam}X
\end{align}
for 
 $X\in \Catlam$
(See \eqref{item:endo-iso}
of 
Theorem \ref{Th:end-ring}).
Then
$G(X)$ is an object of $\M^{\slam} $. 
Here the action of $D$ on $G(X)$
is defined in an obvious way.


\begin{Th}\label{Th:equivalent-cat}
Let $\lam\in \Paf{\cri}$.
\begin{enumerate}
 \item \label{item:free}
The $V(\lam)$ is a free $\RZlam$-module.
\item {\rm{(}}cf.\ Theorem 6.3 of \cite{FreGai04}{\rm{)}}
  The functor $F$
gives an equivalence 
of categories 
$\M^{\slam}\isomap  \Catlam$.
The inverse functor is given by
$G$.
\end{enumerate}
\end{Th}
 The proof of
Theorem \ref{Th:equivalent-cat} 
is given in \cite{Ara08}.
%

The following assertion 
follows immediately
from 
 Theorem \ref{Th:equivalent-cat}.
\begin{Th}\label{Th:chonj-tru}
 Conjecture \ref{Conj-FG} holds.
\end{Th}

Because 
$L(\lam)=G(\C_{\chilam})$,
by \eqref{item:free} of Theorem \ref{Th:equivalent-cat}
one has
\begin{align}
 \ch V(\lam)=\ch \RZlam\cdot  \ch L(\lam).
\end{align}
Therefore 
\eqref{eq:2007-05-28-00-08}
and \eqref{item:ch-of-end} of Theorem \ref{Th:end-ring}
give
Theorem \ref{Th:Main}.


\begin{Rem}
Using  Theorem \ref{Th:equivalent-cat},
one can  show 
the irreduciblity of 
the $\g'$-module
$V(\lam)/\ker \chi \cdot  V(\lam)$ 
for any $\chi\in \Op_{{}^LG}^{\slam}$.
This confirms the original conjecture of
Frenkel and Gaitsgory (Conjecture 5 of \cite{Fre06},
see Remarks \ref{Rem:originial-conj}
and \ref{Rem:reg-op}) partially.
\end{Rem}\subsection{Acknowledgment}
This work grew out of stimulating conversations
with E. Frenkel, 
and
the author is grateful to him.
He also thanks for W.  Wang for 
useful discussions.

\bibliographystyle{alpha}
\bibliography{math}
\end{document}